\documentclass[10pt]{article}
\usepackage{amsmath,amsthm,amssymb}
\usepackage{graphicx}

\title{Ramanujan's cubic transformation inequalities for zero-balanced hypergeometric functions}
\author{\small Miao-Kun Wang$^{a}$, Yu-Ming Chu$^{a}$, Ye-Ping Jiang$^{b}$}
\date{}

\begin{document}
\maketitle
\renewcommand{\thefootnote}{\fnsymbol{footnote}}

{\footnotesize\rm \noindent $^a$ Department of Mathematics, Huzhou Teachers College, Huzhou
313000, China;\\
$^b$ College of Mathematics and Econometrics, Hunan University, Changsha 410082, China.\\
Correspondence should be addressed to Yu-Ming Chu, chuyuming@hutc.zj.cn}

\medskip
\noindent{\bf Abstract}:  In this paper, a generalization of Ramanujan's cubic transformation, in the form of an inequality, is proved for zero-balanced Gaussian hypergeometric function $F(a,b;a+b;x)$, $a,b>0$.

\noindent{\bf Keywords}: Gaussian hypergeometric function,
Ramanujan's cubic transformation, inequality

\noindent{\bf 2010 Mathematics Subject Classification}: 33C05

\bigskip
\noindent{\bf\large 1. Introduction}
\bigskip
\setcounter{section}{1}
\setcounter{equation}{0}

For real numbers $a$, $b$ and $c$ with $c\neq 0, -1, -2, \cdots$,
the Gaussian hypergeometric function is defined by
\begin{equation}
F(a,b;c;x)={}_{2}F_{1}(a,b;c;x)=\sum_{n=0}^{\infty}\frac{(a,n)(b,n)}{(c,n)}\frac{x^n}{n!}
\end{equation}
for $x\in(-1,1)$, where $(a,n)$ denotes the shifted factorial
function $(a,n)=a(a+1)(a+2)(a+3)\cdots(a+n-1)$ for $n=1,2,\cdots$,
and $(a,0)=1$ for $a\neq 0$. And $F(a,b;c;x)$ is called
zero-balanced if $c=a+b$.

It is well known that $F(a,b;c;x)$ has many important applications
in various fields of the mathematical and natural sciences [1-2],
and many classes of special function in mathematical physics are
particular cases of this function [3]. For a extensive list of
$F(a,b;c;x)$ see [4-7].

As the special case of Gaussian hypergeometric
function, for $r\in(0,1)$, Legendre's complete elliptic integrals of the first kind
is defined by
\begin{equation*}
{\mathcal{K}}(r)=\int_{0}^{{\pi}/{2}}(1-r^2\sin^2\theta)^{-1/2}d
\theta=\frac{\pi}{2}F(\frac{1}{2},\frac{1}{2};1;r^2).
\end{equation*}
Some of the most important properties of the elliptic integrals $\mathcal{K}(r)$ are the Landen identities:
\begin{equation*}
\mathcal{K}\left(\frac{2\sqrt{r}}{1+r}\right)=(1+r)\mathcal{K}(r),\quad
\mathcal{K}\left(\frac{1-r}{1+r}\right)=\frac{1+r}{2}\mathcal{K}(\sqrt{1-r^2}),
\end{equation*}
namely,
\begin{equation}
F(\frac{1}{2},\frac{1}{2};1;\frac{4r}{(1+r)^2})=(1+r)F(\frac{1}{2},\frac{1}{2};1;r^2),
\end{equation}
\begin{equation}
F(\frac{1}{2},\frac{1}{2};1;\left(\frac{1-r}{1+r}\right)^2)=\frac{1+r}{2}F(\frac{1}{2},\frac{1}{2};1;1-r^2).
\end{equation}

For zero-balanced Gaussian hypergeometric functions $F(a,b;a+b;x)$,
$a,b>0$, Simi\'{c} and Vuorinen [8] determined maximal region of
$ab$ plane where equations (1.2) and (1.3) turn on respective
inequalities valid for each $x\in(0,1)$.

As is known to all, Ramanujan's cubic transformation is defined as
\begin{equation}
F(\frac{1}{3},\frac{2}{3};1;1-\left(\frac{1-r}{1+2r}\right)^3)=(1+2r)F(\frac{1}{3},\frac{2}{3};1;r^3),
\end{equation}
\begin{equation}
F(\frac{1}{3},\frac{2}{3};1;\left(\frac{1-r}{1+2r}\right)^3)=\frac{1+2r}{3}F(\frac{1}{3},\frac{2}{3};1;1-r^3).
\end{equation}

Inspired by the ideas of Simi\'{c} and Vuorinen [8], we find the
maximal region of $ab$ plane for $F(a,b;a+b;x)$, $a,b>0$ where
equations (1.4) and (1.5) turn on respective inequalities valid for
each $x\in(0,1)$.

The following asymptotic formulas for zero-balanced hypergeometric function (see [9, 10]) will be used in this
paper.
\begin{equation}
F(a,b;a+b;r)\sim -\frac{1}{B(a,b)}\log(1-r)
\end{equation}
and
\begin{equation}
B(a,b)F(a,b;a+b;r)+\log(1-r)=R(a,b)+O((1-r)\log(1-r)),
\end{equation}
as $r$ tends to 1, where
\begin{equation}
B(z,w)=\frac{\Gamma(z)\Gamma(w)}{\Gamma(z+w)},\quad \mbox{Re}\
{z}>0,\quad \mbox{Re}\ {w}>0
\end{equation}
is the classical beta function,
\begin{equation}
R(a,b)=-\Psi(a)-\Psi(b)-2\gamma,\quad R(1/3,2/3)=\log{27},
\end{equation}
\begin{equation}
\Psi(z)=\frac{d}{dz}(\log\Gamma(z))=\frac{\Gamma'(z)}{\Gamma(z)},\quad
\mbox{Re}\ {z}>0,
\end{equation}
and $\gamma$ is the Euler-Mascheroni constant.

\medskip
{\bf Lemma 1.1} (See [8, Lemma 1.1]). Suppose that the power series
$f(x)=\sum\limits_{n=0}^{\infty}a_{n}x^{n}$ and
$g(x)=\sum\limits_{n=0}^{\infty}b_{n}x^{n}$ have the radius of
convergence $r>0$ and $b_{n}>0$ for all $n\in\{0,1,2,\cdots\}$. Let
$h(x)={f(x)}/{g(x)}$, then

(1) If the sequence $\{a_{n}/b_{n}\}_{n=0}^{\infty}$ is (strictly)
increasing (decreasing), then $h(x)$ is also (strictly) increasing
(decreasing) on $(0,r)$;

(2) If the sequence $\{a_{n}/b_{n}\}$ is (strictly) increasing
(decreasing) for $0<n\leq n_{0}$ and (strictly) decreasing
(increasing) for $n>n_{0}$, then there exists $x_{0}\in(0,r)$ such
that $h(x)$ is (strictly) increasing (decreasing) on $(0,x_{0})$ and
(strictly) decreasing (increasing) on $(x_{0},r)$.

\bigskip
\noindent{\bf\large 2. Main Results}
\bigskip
\setcounter{section}{2} \setcounter{equation}{0}

For convenience, we first introduce the following regions
in $\{(a,b)\in\mathbb{R}^2|a>0,b>0\}$:

\begin{equation*}
D_{1}=\{(a,b)|a,b>0, ab\leq 2/9, ab-\frac{2}{9}(a+b)\leq 0\},
\end{equation*}
\begin{equation*}
D_{2}=\{(a,b)|a,b>0, ab<2/9, ab-\frac{2}{9}(a+b)> 0\},
\end{equation*}
\begin{equation*}
D_{3}=\{(a,b)|a,b>0, ab\geq 2/9, ab-\frac{2}{9}(a+b)\geq 0\},
\end{equation*}
\begin{equation*}
D_{4}=\{(a,b)|a,b>0, ab>2/9, ab-\frac{2}{9}(a+b)<0\},
\end{equation*}
\begin{equation*}
D_{5}=\{(a,b)|a,b>0, a+b\leq 1, ab-\frac{2}{9}(a+b)\leq 0\},
\end{equation*}
\begin{equation*}
D_{6}=\{(a,b)|a,b>0, a+b\geq 1, ab-\frac{2}{9}(a+b)\geq 0\}.
\end{equation*}
Clearly, $D_{1}\cup D_{2}\cup D_{3}\cup D_{4}=\{(a,b)\in\mathbb{R}^2|a>0,b>0\}$, $D_{5}\subset D_{1}$ and $D_{6}\subset D_{3}$.

\medskip
{\bf Theorem 2.1.} If $(a,b)\in D_{1}$, then the inequality
\begin{equation}
F(a,b;a+b;\frac{9r(1+r+r^2)}{(1+2r)^3})\leq (1+2r)F(a,b;a+b;r^3)
\end{equation}
holds for all $r\in(0,1)$. Also, if $(a,b)\in D_{3}$, then the
reversed inequality
\begin{equation}
F(a,b;a+b;\frac{9r(1+r+r^2)}{(1+2r)^3})\geq (1+2r)F(a,b;a+b;r^3)
\end{equation}
takes place for each $r\in (0,1)$, with equality in each instance if and only if $(a,b)=(1/3,2/3)$ or $(a,b)=(2/3,1/3)$.

In the remaining region $(a,b)\in D_{2}\cup D_{4}$, neither of the
above inequalities holds for each $r\in(0,1)$.

\medskip
{\bf Theorem 2.2.} If $(a,b)\in D_{1}$, then the double inequality
\begin{equation}
1\leq
\frac{(1+2r)F(a,b;a+b;r^3)}{F(a,b;a+b;\frac{9r(1+r+r^2)}{(1+2r)^3})}\leq
\frac{\sqrt{3}B(a,b)}{2\pi}
\end{equation}
holds for all $r\in(0,1)$. And if $(a,b)\in D_{3}$, then inequality (2.3) is reversed
\begin{equation}
\frac{\sqrt{3}B(a,b)}{2\pi}\leq
\frac{(1+2r)F(a,b;a+b;r^3)}{F(a,b;a+b;\frac{9r(1+r+r^2)}{(1+2r)^3})}\leq
1.
\end{equation}
Moreover, both bounds in inequalities (2.3) and (2.4) are sharp and
each equality is reached for $a=1/3$ and $b=2/3$, or $a=2/3$ and
$b=1/3$.

\medskip
{\bf Corollary 2.3.} For $r\in(0,1)$, and $(a,b)\in D_{1}$, one has
\begin{equation}
\frac{2\pi}{\sqrt{3}}\frac{1}{B(a,b)}F(a,b;a+b;r^3)<F(a,b;a+b;\frac{9r(1+r+r^2)}{(1+2r)^3})<3F(a,b;a+b;r^3).
\end{equation}
In the region $(a,b)\in D_{3}$, one has
\begin{equation}
F(a,b;a+b;r^3)<F(a,b;a+b;\frac{9r(1+r+r^2)}{(1+2r)^3})<\frac{6\pi}{\sqrt{3}}\frac{1}{B(a,b)}F(a,b;a+b;r^3).
\end{equation}

\medskip
{\bf Theorem 2.4.} Let $B=B(a,b)$ and $R=R(a,b)$ are defined as in
(1.8) and (1.9), respectively. Then for $(a,b)\in D_{5}$, inequality
\begin{equation}
0\leq
(1+2r)F(a,b;a+b;r^3)-F(a,b;a+b;\frac{9r(1+r+r^2)}{(1+2r)^3})\leq
\frac{2(R-\log{27})}{B}
\end{equation}
holds for all $r\in(0,1)$. Also, for $(a,b)\in D_{6}$,
\begin{equation}
0\leq
F(a,b;a+b;\frac{9r(1+r+r^2)}{(1+2r)^3})-(1+2r)F(a,b;a+b;r^3)\leq
\frac{2(R-\log{27})}{B}.
\end{equation}

\medskip
{\bf Theorem 2.5.} For $(a,b)\in D_{1}$ and each $x\in(0,1)$, one
has
\begin{equation}
\frac{1}{3}\leq \frac{F(a,b;a+b;\left(\frac{1-x}{1+2x}\right)^3)}{(1+2x)F(a,b;a+b;1-x^3)}\leq \frac{\sqrt{3}B(a,b)}{6\pi}.
\end{equation}

(2) For $(a,b)\in D_{3}$ and each $x\in(0,1)$, one has
\begin{equation}
\frac{\sqrt{3}B(a,b)}{6\pi}\leq \frac{F(a,b;a+b;\left(\frac{1-x}{1+2x}\right)^3)}{(1+2x)F(a,b;a+b;1-x^3)}\leq \frac{1}{3}.
\end{equation}

(3) For $(a,b)\in D_{5}$ and each $x\in(0,1)$, we have
\begin{align}
(1+2x)F(a,b;a+b;1-x^3)\leq
3F(a,b;a+b;\left(\frac{1-x}{1+2x}\right)^3)\nonumber\\
\leq (1+2x)\left[F(a,b;a+b;1-x^3)+\frac{2(R(a,b)-\log{27})}{B(a,b)}\right].
\end{align}

(4) For $(a,b)\in D_{6}$ and each $x\in(0,1)$, we have
\begin{align}
&0\leq(1+2x)F(a,b;a+b;1-x^3)-
3F(a,b;a+b;\left(\frac{1-x}{1+2x}\right)^3)\nonumber\\
&\qquad\qquad\qquad\leq \frac{2(1+2x)(\log{27}-R(a,b))}{B(a,b)}.
\end{align}

\bigskip
\noindent{\bf\large 3. Proofs of Theorems}
\bigskip
\setcounter{section}{3} \setcounter{equation}{0}

In order to prove our main results, we introduce several symbols.
Throughout this section, we let
\begin{equation*}
F(x)=F(a,b;a+b;x),\quad G(x)=F(a,b;a+b+1;x),
\end{equation*}
where $a,b>0$ with $(a,b)\neq (1/3,2/3)$ and $(a,b)\neq (2/3,1/3)$, and
\begin{equation*}
F^{*}(x)=F(\frac{1}{3},\frac{2}{3};1;x),\quad G^{*}(x)=F(\frac{1}{3},\frac{2}{3};2;x).
\end{equation*}

\medskip
{\bf Lemma 3.1.} (1) The function $f(r)=F(r)/F^*(r)$ is strictly decreasing in $(0,1)$ on $D_{1}$, and strictly increasing in $(0,1)$ on $D_{3}$. Moreover, if $(a,b)\in D_{2}$ ($D_{4}$, resp.), then there exists $r_{0}$ ($r_{0}^*$, resp.) such that $f(r)$ is strictly increasing (decreasing, resp.) in $(0,r_{0})$ ($(0,r_{0}^*)$, resp.), and strictly decreasing (increasing, resp.) in $(r_{0},1)$ ($(r_{0}^*,1)$, resp.);

(2) The function $g(r)=G(r)/G^{*}(r)$ is strictly decreasing in
$(0,1)$ on $D_{5}$ and strictly increasing in $(0,1)$ on $D_{6}$.

\medskip
{\bf Proof.} For part (1), denote by $A_{n}=(a,n)(b,n)/[(a+b,n)n!]$
and $A_{n}^*=(1/3,n)(2/3,n)/[(n)!]^2$, then
\begin{equation}
f(r)=\frac{F(r)}{F^{*}(r)}=\frac{\sum\limits_{n=0}^{\infty}A_{n}r^n}{\sum\limits_{n=0}A_{n}^*r^n}.
\end{equation}
Note that the monotonicity of $\{A_{n}/A_{n}^*\}$ depends on the sign of
\begin{equation}
H_{n}=(ab-\frac{2}{9})n+ab-\frac{2}{9}(a+b).
\end{equation}

We divide the proof into four cases.

{\bf Case 1} $(a,b)\in D_{1}$. Then (3.2) implies $H_{n}<0$ for $n=0,1,2,\cdots$, and $f(r)$ is strictly decreasing in $(0,1)$ by (3.1) and Lemma 1.1.

{\bf Case 2} $(a,b)\in D_{3}$. Then (3.2) implies $H_{n}>0$ for $n=0,1,2,\cdots$, and $f(r)$ is strictly increasing in $(0,1)$ by (3.1) and Lemma 1.1.

{\bf Case 3} $(a,b)\in D_{2}$. Then from (3.2) we conclude that the sequence $\{A_{n}/A_{n}^*\}$ increases and then decreases. By (3.1) and  Lemma 1.1(3), there exists $r_{0}\in(0,1)$ such that $f(r)$ is strictly increasing in $(0,r_{0})$ and strictly decreasing in $(r_{0},1)$.

{\bf Case 4} $(a,b)\in D_{4}$. Then from (3.2) we know that the sequence $\{A_{n}/A_{n}^*\}$ decreases and then increases. By (3.1) and Lemma 1.1(3), there exists $r_{0}^*\in(0,1)$ such that $f(r)$ is strictly decreasing in $(0,r_{0}^*)$ and strictly increasing in $(r_{0}^*,1)$.

For part (2), denote by $B_{n}=(a,n)(b,n)/[(a+b+1,n)n!]$ and
$B_{n}^*=(1/3,n)(2/3,n)/[(2,n)(n)!]$, then
\begin{equation}
g(r)=\frac{G(r)}{G^{*}(r)}=\frac{\sum\limits_{n=0}^{\infty}B_{n}r^n}{\sum\limits_{n=0}B_{n}^*r^n}.
\end{equation}
Note that the monotonicity of $\{B_{n}/B_{n}^*\}$ depends on the sign of
\begin{equation}
H_{n}^*=(a+b+ab-\frac{11}{9})n+\frac{2}{9}(9ab-a-b-1).
\end{equation}

We divide the proof into two cases.

{\bf Case A} $(a,b)\in D_{5}$. Then
$a+b+ab-11/9\leq 11(a+b)/9-11/9\leq 0$ and  $9ab-a-b-1=9ab-2(a+b)+(a+b)-1\leq 0$.
Thus $H_{n}^*<0$ for $n=0,1,2,\cdots$ (because $(a,b)\neq (1/3,2/3)$ and $(a,b)\neq (2/3,1/3)$) by (3.4). Therefore, $g(r)$ is strictly decreasing in $(0,1)$ follows from (3.3) and Lemma 1.1.

{\bf Case B} $(a,b)\in D_{6}$. Then $a+b+ab-11/9\geq 11(a+b)/9-11/9\geq 0$ and $9ab-a-b-1=9ab-2(a+b)+(a+b)-1\geq 0$.
Thus $H_{n}^*>0$ for $n=0,1,2,\cdots$ by (3.4). Therefore, $g(r)$ is strictly increasing in $(0,1)$ follows from (3.3) and Lemma 1.1.

\medskip
{\bf Proof of Theorem 2.1.} Let $x=x(r)=r^3$ and $y=y(r)=9r(1+r+r^2)/(1+2r)^3$, then simple computation leads to $0<x<y<1$ for $0<r<1$. Using Lemma 3.1(1), we get $f(x)>f(y)$ on $D_{1}$, and $f(x)<f(y)$ on $D_{3}$.

For $(a,b)\in D_{1}$, by (1.4), one has
\begin{equation*}\
\frac{F(r^3)}{F^*(r^3)}>\frac{F(y)}{F^*(y)},\quad F(y)<\frac{F^{*}(y)}{F^*(r^3)}F(r^3)=(1+2r)F(r^3).
\end{equation*}
Thus inequality (2.1) follows.

Inequality (2.2) is obtained analogously. The remaining conclusions easily follows from Lemma 3.1(1).

\medskip
{\bf Proof of Theorem 2.2.}  Let $f(r)$ be defined as in Lemma 3.1(1), then $f(r)$ is strictly decreasing on $D_{1}$. Then (1.6) leads to
\begin{equation*}
1=\lim_{r\rightarrow 0^+}\frac{F(r)}{F^*(r)}>\frac{F(r)}{F^*(r)}>\lim_{r\rightarrow 1^-}\frac{F(r)}{F^*(r)}=\frac{B(1/3,2/3)}{B(a,b)}=\frac{2\sqrt{3}\pi}{3B(a,b)}
\end{equation*}
and
\begin{equation*}
\frac{\sqrt{3}B(a,b)}{2\pi}\frac{1}{F^*(y(r))}>\frac{1}{F(y(r))}\Longrightarrow \frac{\sqrt{3}B(a,b)}{2\pi}\frac{F^*(x(r))}{F^*(y(r))}>\frac{F(x(r))}{F(y(r))}.
\end{equation*}
Thus inequality (2.3) is clear.

Inequality (2.4) valid on $D_{3}$ can be proved similarly.

\medskip
{\bf Lemma 3.2.} The function
\begin{equation*}
J(r)=(1+2r^{1/3})F(a,b;a+b;r)-F(a,b;a+b;\frac{9r^{1/3}(1+r^{1/3}+r^{2/3})}{(1+2r^{1/3})^3})
\end{equation*}
is strictly increasing in $(0,1)$ on $D_{5}$ and strictly decreasing in $(0,1)$ on $D_{6}$.

\medskip
{\bf Proof.} Let $z=9r^{1/3}(1+r^{1/3}+r^{2/3})/(1+2r^{1/3})^3$. Then
\begin{equation*}
1-z=\frac{(1-r^{1/3})^3}{(1+2r^{1/3})^3},\quad \frac{dz}{dr}=\frac{3(1-r^{1/3})^2}{r^{2/3}(1+2r^{1/3})^4}.
\end{equation*}
Note that
\begin{equation*}
(1-x)F(a+1,b+1;a+b+1;x)=F(a,b;a+b+1;x).
\end{equation*}

Differentiating $J(r)$ gives
\begin{align}
r^{2/3}(1-r^{1/3})J'(r)=&\frac{2}{3}(1-r^{1/3})F(a,b;a+b;r)+\frac{ab}{a+b}\frac{r^{2/3}(1+2r^{1/3})(1-r^{2/3})}{1-r}\nonumber\\
&\times F(a,b;a+b+1;r)-\frac{3ab}{(a+b)(1+2r^{1/3})}F(a,b;a+b+1;z)\nonumber\\
=&\frac{2}{3}(1-r^{1/3})F(r)+\frac{ab}{a+b}\frac{r^{2/3}(1+2r^{1/3})(1-r^{2/3})}{1-r}G(r)\nonumber\\
&-\frac{3ab}{(a+b)(1+2r^{1/3})}G(z).
\end{align}

On the other hand, differentiating Ramanujan cubic transformation, we get
\begin{equation}
\frac{2}{3}\frac{G^{*}(z)}{1+2r^{1/3}}=\frac{2}{3}(1-r^{1/3})F^{*}(r)+\frac{2}{9}\frac{r^{2/3}(1+2r^{1/3})(1-r^{2/3})}{1-r}G^*(r).
\end{equation}

Let $g(r)$ be defined as in Lemma 3.1(2), then $g(r)$ is strictly decreasing in $(0,1)$ on $D_{5}$. Then from $0<r<z<1$ we get $g(r)>g(z)$, namely
\begin{equation}
G(z)<\frac{G^{*}(z)}{G^*(r)}G(r).
\end{equation}
Equations (3.5) and (3.6) together with inequality (3.7) yield
\begin{align*}
&r^{2/3}(1-r^{1/3})J'(r)>\frac{2}{3}(1-r^{1/3})F(r)+\frac{ab}{a+b}\frac{r^{2/3}(1+2r^{1/3})(1-r^{2/3})}{1-r}G(r)\\
&-\frac{3ab}{(a+b)(1+2r^{1/3})}\frac{G^*(z)}{G^*(r)}G(r)\\
&=\frac{2}{3}(1-r^{1/3})F(r)+\frac{ab}{a+b}\frac{r^{2/3}(1+2r^{1/3})(1-r^{2/3})}{1-r}G(r)-\frac{3ab}{(a+b)(1+2r^{1/3})}\\
&\times\left((1-r^{1/3})(1+2r^{1/3})\frac{F^*(r)}{G^*(r)}+\frac{1}{3}\frac{r^{2/3}(1+2r^{1/3})^2(1-r^{2/3})}{1-r}\right)G(r)\\
&=\frac{2}{3}(1-r^{1/3})F(r)-\frac{3ab}{(a+b)}(1-r^{1/3})\frac{F^*(r)}{G^*(r)}G(r)\\
&=\frac{2}{3}(1-r^{1/3})\left[F(r)-\frac{9ab}{2(a+b)}\frac{F^*(r)}{G^*(r)}G(r)\right].
\end{align*}
Note that
\begin{equation*}
\frac{F'(r)}{{F^{*}}'(r)}=\frac{9ab}{2(a+b)}\frac{G(r)}{G^{*}(r)}.
\end{equation*}
Thus
\begin{equation*}
\frac{3}{2}r^{2/3}J'(r)>F(r)-\frac{F'(r)}{{F^{*}}'(r)}F^{*}(r)=\frac{F^{2}(r)}{{F^{*}}'(r)}\left(\frac{F^{*}(r)}{F(r)}\right)'.
\end{equation*}

It follows from Lemma 3.1(1) and the fact that $D_{5}\subset D_{1}$ that $\left({F^{*}(r)}/{F(r)}\right)'\geq 0$ on $D_{5}$. Hence
$J'(r)>0$, and $J(r)$ is strictly increasing in $(0,1)$ on $D_{5}$.

Since $g(r)$ is strictly decreasing in $(0,1)$ on $D_{6}$, we have $g(r)>g(z)$, namely
\begin{equation*}
G(z)>\frac{G^*(z)}{G^*(r)}G(r).
\end{equation*}
With the similar argument, one has
\begin{equation*}
\frac{3}{2}r^{2/3}J'(r)<\frac{F^{2}(r)}{{F^{*}}'(r)}\left(\frac{F^{*}(r)}{F(r)}\right)'<0,
\end{equation*}
since $f(r)=F(r)/F^*(r)$ is strictly increasing in $(0,1)$ on $D_{6}\subset D_{3}$. Hence $J(r)$ is strictly decreasing in $(0,1)$ on $D_{6}$. $\Box$

\medskip
{\bf Proof of Theorem 2.4.} By Lemma 3.2 we obtain $\lim\limits_{r\rightarrow 0^+}J(r)<J(r)<\lim\limits_{r\rightarrow 1^-}J(r)$ on $D_{5}$ and $\lim\limits_{r\rightarrow 1^-}J(r)<J(r)<\lim\limits_{r\rightarrow 0^+}J(r)$ on $D_{6}$. Clearly, $\lim\limits_{r\rightarrow 0^+}J(r)=0$. And by (1.7), we have
\begin{align*}
&\lim\limits_{r\rightarrow 1^-}J(r)\\
=&\lim\limits_{r\rightarrow 1^-}\frac{3R(a,b)-3\log(1-r)-(R(a,b)-3\log[(1-r^{1/3})/(1+2r^{1/3})])+o(1)}{B(a,b)}\\
=&\frac{2(R(a,b)-\log{27})}{B(a,b)}.
\end{align*}

The assertion of Theorem 2.4 follows.

\medskip
{\bf Proof of Theorem 2.5.}  Theorem 2.5 follows from Theorems 2.2 and 2.4 with $x=(1-r)/(1+2r)\in(0,1)$.

\medskip
\emph{\textbf{Acknowledgement}.} This research was supported by the Natural Science
Foundation of China (Grant Nos. 11071059, 11071069, 11171307), and
the Innovation Team Foundation of the
Department of Education of Zhejiang Province (Grant no. T200924).

\bigskip
\def\refname
\end{document}